\documentclass[12pt,a4paper]{article}

\usepackage{amssymb,amsmath}
\usepackage[mathscr]{eucal}  
\usepackage{amsfonts,amssymb,latexsym,amsmath,amsthm,epsfig,amscd}

\newtheorem{lemma}{Lemma}[section]
\newtheorem{cor}[lemma]{Corollary}

\newtheorem{theorem}[lemma]{Theorem}
\newtheorem{pro}[lemma]{Proposition}

\newtheorem {defn}{Definition}
\newtheorem{remark}{Remark}

\newcommand{\leps}{{\cal L}{^{^{(\epsilon_0)}}}}
\newcommand{\RN}{\mathbb{R}^N}
\newcommand{\K}{\text{supp} (\mu) }
\newcommand{\N}{N_{X}(\beta ,\x,r) }
\newcommand{\x}{\mathbf{x}}
\newcommand{\y}{\mathbf{y}}

\newcommand{\p}{\mathbf{p}}
\newcommand{\A}{{\cal{A}}}

\begin{document}
\title{Schmidt's Game on Fractals}
\author{Lior Fishman} 
\maketitle

\begin{abstract}

We construct $(\alpha ,\beta)$ and $\alpha $-winning sets 
in the sense of Schmidt's game, played on the support of certain measures 
(absolutely friendly) and show how to compute the Hausdorff dimension for some. 

In particular we prove that if $K$ is the attractor of an irreducible finite family 
of contracting similarity maps of $\RN $ satisfying the open set condition, 
(the Cantor's ternary set, Koch's curve and Sierpinski's gasket to name a few
known examples), then for any countable collection of non-singular affine transformations, 
$\Lambda_i:\RN \rightarrow \RN $, 
\begin{center}
dim$K=$ dim$K\cap (\cap ^{\infty}_{i=1}( \Lambda_i(\textbf{BA})))$
\end{center}
where \textbf{BA} is the set of badly approximable vectors in $\RN$.       
\end{abstract}   
\vspace{2cm} 

\setcounter{section}{-1}

\section{Introduction}

We shall be using Schmidt's game first introduced by W. M. Schmidt \cite{S} for 

\noindent estimating the Hausdorff dimension of certain sets. 
Let us first define the set of badly approximable vectors.
A vector $\x\in\RN $ is said to be badly approximable if there exists 
$\delta >0$ such that for any $\p \in \mathbb{Z}^{N}$, $q\in \mathbb{N}^{+}$ 
\begin{equation}
d(\x,\frac{\p}{q})\geq \delta q^{-\frac{N+1}{N}}
\end{equation}
where $d$ is the Euclidean distance function between points.
We denote the set of all badly approximable vectors by \textbf{BA}.
The above mentioned game was used by Schmidt, among other things, 
to tackle the following questions concerning \textbf{BA}:
\begin{enumerate}
	\item If $\left\{\Lambda_i\right\}^{\infty}_{i=0}$ is a 
countable collection of non-singular affine transformations 
	
	$\Lambda_i:\RN \rightarrow \RN $, is $\cap ^{\infty}_{i=1}( \Lambda_i(\textbf{BA}))\neq \emptyset$?
	\item If $\cap ^{\infty}_{i=1}( \Lambda_i(\textbf{BA}))\neq \emptyset$, what is 
dim$\cap ^{\infty}_{i=1}( \Lambda_i(\textbf{BA}))$?
\end{enumerate}
 
Schmidt proved not only that the intersection in non empty, 
but is in fact ``large'' dimension wise, i.e., is of dimension $N$.

In recent years similar questions have been posed regarding the intersection of \textbf{BA} 
with certain subsets of $\RN $. For example, let $K$ be any of the following sets: 
Cantor's ternary set, Koch's curve, Sierpinski's gasket, 
or in general, an attractor of an irreducible finite family of contracting similarity maps 
of $\RN $ satisfying the open set condition. 
(This condition due to J. E. Hutchinson \cite{H} is discussed in section 5). 
One may ask the following questions:
\begin{enumerate}
	\item Is $K\cap \textbf{BA}\neq \emptyset$?
	\item If $K\cap \textbf{BA}\neq \emptyset$, what is dim$K\cap \textbf{BA}$?
\end{enumerate}
Answers to both of these questions have been independently given in \cite {KW} and \cite{KTV} 
proving dim$K\cap \textbf{BA}$=dim$K$ for a large family of sets including those mentioned above.

This paper's aim is to extend these results, utilizing Schmidt's game, by 

\noindent answering the following question:
If $\left\{\Lambda_i\right\}^{\infty}_{i=0}$ is a countable collection of 
non-singular affine transformations 
$\Lambda_i:\RN \rightarrow \RN $, what is dim$K\cap (\cap ^{\infty}_{i=1}( \Lambda_i(\textbf{BA})))$?

It turns out that for a large family of sets 
the answer is analogous to Schmidt's result in $\RN$, 
namely we prove in section 5, 

\vspace{5 mm}

\textsl{Corollary \ref{MAIN EXAMPLE}.
Let $\left\{ \phi_1,...,\phi_k\right)$ be a finite irreducible family of 
contracting similarity maps of $\RN $ satisfying the open set condition and
let $K$ be its attractor.
Then for any countable collection 
of non-singular affine transformations 
\noindent $\left\{\Lambda_i\right\}^{\infty}_{i=0}$}, 
\textsl{
\noindent with $\Lambda_i:\RN \rightarrow \RN $ the set
\begin{center} 
${\cal{S}}= K\cap (\cap ^{\infty}_{i=1}( \Lambda_i(\textbf{BA})))$ 
\end{center} 
is a winning set on K.} 
\textsl{Furthermore, dim${\cal \cal{S}}$=dim$K$.}

\vspace{5 mm}

Our research closely follows in the footsteps of \cite{KLW}, \cite{KW}
and consequently \cite{PV} and \cite{KTV}.
The definitions of measures given in the first and third of the above
mentioned papers were 
not originally intended for creating a ``friendly'' environment for Schmidt's game 
on their support. 
It turns out however that in a sense to be made clearer later, these measures indeed 
provide an hospitable playground for this game.

\vspace{1.5cm}

Section 1 is devoted to establishing the link between the definitions given in \cite{KLW}, the
stronger assumptions in \cite{PV} and our work,
exhibiting a geometric feature material for later discussion.

\vspace{5 mm}

In section 2 we follow the general setup introduced in \cite{KTV}
proving as a consequence of corollary \ref{BA=BAD} and theorem \ref{BAD=ALPHABETA} 
that if a measure $\mu$ is absolutely friendly (see definition in section 1) then
under certain conditions
\begin{center} 
$BA\cap \K$ is an $(\alpha , \beta)$-winning set on $\K$.
\end{center}

\vspace{5 mm}

In section 3 we formulate a sufficient condition
for establishing a lower bound of a winning set's Hausdorff dimension,
where the winning set is a subset of the support of an absolutely friendly measure.

\vspace{5 mm}

In section 4 we prove an analogue to the simplex lemma in \cite{S}.

\vspace{5 mm}

Section 5 is our main example, an application to the Hutchinson construction. 

\vspace{5 mm}

\noindent As should be obvious from the discussion above, our conclusions strengthen
results in \cite{KW} and \cite{KTV} regarding the Hausdorff dimension of the 
intersection of BA with certain sets. (See corollary 1.2 in \cite{KW}
and conclusions from theorem 1 in \cite{KTV}). 
We should however note that in proving our theorems we are in fact using
stronger assumptions on our measures in order to make sure that our target set-the set
of badly approximable vectors, is indeed a winning set on the support of these measures.

\vspace{2cm}

\subsection*{Notation}

\noindent $\mathbb{R}$, $\mathbb{Q}$ and $\mathbb{N}$ denote the set of real, 
rational and natural numbers respectively.\\

\noindent $\mathbb{R}^{+}$ is the set of non-negative 
real numbers while $\mathbb{N}^{+}$ denotes the set of strictly positive integers. \\

\noindent Boldface lower case letters ($\x$, $\y$,...etc.) denote points in $\mathbb{R}^N$.\\

\noindent The function $d$ is the Euclidean distance function between points. 
If $A$ and $B$ are any two subsets of $\mathbb{R}^N$, $d(A,B)
= \text{inf}\left\{d(\x,\y):\x\in A,\y\in B\right\}$.\\

\noindent $\lambda _N$ denotes the Lebesgue measure in $\RN $.\\

\noindent In the metric space $(\RN ,d)$,  
$B(x,r)$ will denote a 
closed ball of radius $r$ centered at $x$, 
i.e., $B(x,r)=\left\{z:d(x,z)\leq r\right\}$, $\partial B(x,r)$ 
the boundary of $B(x,r)$, i.e., $\left\{z:d(z,x)=r\right\}$ 
and $int B(x,r)$ denotes the interior of $B(x,r)$ i.e., $\left\{z:d(x,z)< r\right\}$.\\

\noindent An affine hyperplane of  $\mathbb{R}^N$ will be denoted by 
$\cal{L}$ while $\cal{L}^{(\epsilon)}$ is defined to be the $\epsilon$ neighborhood of 
$\cal{L}$, i.e., ${\cal{L}}^{(\epsilon)}= \{\x \in\mathbb{R}^N :d(\x,{\cal{L}})\leq\epsilon\}$ 
where $\epsilon$ is a non-negative, possibly zero, real number.\\

\noindent Unless otherwise stated, constants are real, strictly positive numbers.\\

\noindent Throughout the paper, $\mu$ will denote a Borel, locally finite measure on $\RN$.\\

\noindent Whenever discussing a measure we denote its support by $\K $.\\

\noindent In order to avoid unnecessary repetitions, all affine transformations 
referred to in this paper are assumed to be non-singular.\\ 

\noindent Following conventional notation, for every $U\subset \RN $ let 

\begin{center}$\left|U\right|=\text{sup}\left\{d(\x,\y):x,y\in U\right\}$.
\end{center}
If $F\subset \RN$, $\delta >0$ and $\left\{U_i\right\}$ is a countable or 
finite collection of sets we say that $\left\{U_i\right\}$ is a $\delta $-cover of $F$ if 

\begin{center}$F\subset \bigcup^{\infty}_{i=1}U_i$\; and 
for every $i$\; $0\leq \left|U_i\right|\leq \delta $.
\end{center}

\noindent If $F\subset \RN $ and $s\geq 0$ then for every $\delta >0$ we define
\begin{center}$H^{s}_{\delta }(F)= \text{inf}
\left\{\sum ^{\infty}_{i=1}\left|U_i\right|^{s}:\left\{U_i\right\}\text{is a}
\; \delta \text{-cover of F} \right\}$
\end{center} 
and
\begin{center}
$H^{s}(F)=\text{lim}_{\delta \rightarrow 0}H^{s}_{\delta }(F)$
\end{center}

\noindent is the $s$-Hausdorff measure.\\

\noindent The Hausdorff dimension of a set $F\subset \RN $ is defined by
\begin{center} $\text{dim}F=\text{inf}\left\{s:H^{s}
(F)=0\right\}=\text{sup}\left\{s:H^{s}(F)=\infty \right\}$.
\end{center}

\vspace{2cm}

\subsection*{Acknowledgments}
\noindent This paper is a part of the author's doctoral dissertation at Ben Gurion university. 
The support of the Israel Science Foundation is gratefully acknowledged.    

\vspace{5mm}

\noindent My deepest thanks to Dmitry Kleinbock for 
carefully reading the paper and providing many helpful suggestions.

\vspace{5mm}

\noindent Last, but in no way least, 
it gives me great pleasure to thank my supervisor Barak Weiss for his countless valuable suggestions and 
for always being absolutely friendly giving them!\\

\section{Absolutely friendly measures}

We first define absolutely friendly measures and
show how it provides the right setting for our work and results.
The class of friendly measures was first introduced in \cite{KLW}, 
followed by the more restrictive $\alpha$-absolutely friendly measures in \cite{PV}. 
The definition of absolutely friendly coincides with that of 
$\alpha$-absolutely friendly, but as the constant $\alpha$ does not seem to 
have any special status in 
any of the formulas we use, we decided to use the term
absolutely friendly instead. 

\begin{defn} \label{def:aad} 
Call a measure $\mu$ on $\mathbb{R}^{N}$ \textbf{absolutely friendly} if the 
following conditions are satisfied:

There exist constants $r_0$, $C$, $D$ and $a$ such that for every $0<r\leq r_0$ and for every $\x \in \K$:
\begin{itemize}
\item[(i)] for any  $0\leq\epsilon\leq r$, and any affine hyperplane $\cal L$,

$\mu(B(\x,r)\cap{\cal L}^{(\epsilon)})< C(\frac{\epsilon}{r})^{^a}\mu (B(\x ,r))$. 

\item[(ii)] $\mu(B(\x,\frac{5}{6} r))> D\mu(B(\x ,r))$.

\end{itemize}

\end{defn}

\noindent Two remarks are in order.

\begin{remark}
Notice that part (ii) of the above definition is equivalent (up to a change of the constant $D$) 
to the so called ``Federer doubling property'' with 
$\frac{1}{2}$ replacing $\frac{5}{6}$. 
\end{remark}

\begin{remark}
The reader should compare (i) with the following more general definition (2.5 in
\cite{KLW}), namely  
given $C$, $a >0$ and an open subset $U$ of $\RN$ we say that 
$\mu $ is 
\textbf{absolutely} $\boldsymbol{(C,a)}$ \textbf{-decaying on U} if for any non-empty open ball $B\subset U$ centered in $\K$, any affine hyperplane ${\cal{L}}\subset \RN$ and any $\epsilon >0$ one has
\begin{equation}
\mu (B\cap {\cal L^{(\epsilon)}})\leq C\left(\frac{\epsilon}{r}\right)^{a}\mu (B)
\end{equation}
where $r$ is the radius of $B$.

\end{remark}
As a consequence of definition \ref{def:aad} we prove the following lemma.

\begin{lemma}\label{GEO}
 
\noindent Suppose $\mu $ is absolutely friendly with constants as in definition \ref{def:aad}.
Define $(\frac{D}{C})^{^{\frac{1}{a}}}= \alpha^{'}$ and let $\cal L$ be any affine hyperplane. 
Then for every $0<r\leq r_0$, if  \;$0<\alpha < \frac{1}{12}\alpha ^{'}$ and 
$0\leq \epsilon_0 <\frac{1}{12}\alpha^{'}r$, 
we have that for every $\x \in \K$ there exists $\x_0\in \K $
such that 
\begin{enumerate}
\item
$B(\x_0,\alpha r)\subset B(\x,r)$
\item
$d(B(\x_0,\alpha r),{\cal L}{^{^{(\epsilon_0)}}})>\alpha r$.
\item
$d(B(\x_0,\alpha r),\partial B(\x,r))>\alpha r$
\end{enumerate}
\end{lemma}

\begin{proof}

\noindent If $d(\x,\leps )>2\alpha r$ the first two conditions 
are evidently satisfied by choosing $\x_0=\x$ while for the third 
notice that $r-\alpha r>\frac{11}{12}r>2\alpha r$.

\noindent Otherwise let $d(\x,\leps )\leq 2\alpha r$.

\noindent Let $\delta =1-\alpha$, \; $\epsilon =5\alpha r +2\epsilon_0$ and denote by 
${\cal L}_\x$ an affine hyperplane parallel to $\cal L$ 
passing through $\x$. We observe that

\begin{equation}
\delta r-\epsilon =(1-6\alpha )r-2\epsilon_0>\left(1-\frac{5}{6}
\alpha^{'}\right)r-\frac{1}{6}\alpha^{'}r=(1-\alpha^{'})r\geq 0
\end{equation}

\begin{equation}
\mu(B(\x,\delta r))=\mu (B(\x,(1-\alpha)r))\geq \mu \left(B\left(\x,\frac {5}{6} 
r\right)\right)\geq D\mu (B(\x,r))
\end{equation}
 
\begin{equation}
\mu ({\cal L}^{(\epsilon)}_\x\cap B(\x,r))\leq C\left(\frac{\epsilon }{r}\right)^{^a}
\mu (B(\x,r))=C\left(5\alpha + \frac{2\epsilon_0}{r}\right) ^{a}\mu (B(\x,r))
\end{equation}

\hspace{1.5cm}$<C(\frac{31}{36}\alpha ^{'})^a\mu (B(\x,r))<C(\alpha ^{'})^a
\mu (B(\x,r))\leq D\mu (B(\x,r)).$\\
 
Consequently, denoting by $\Xi = B(\x,\delta r)-{\cal L}^{(\epsilon)}_\x$, we have 
$\mu (\Xi \cap B(\x,r))>0$ 
and we may choose $\x_0$ to be any point in $\Xi \cap \K$.

The first condition is fulfilled by our choice of $\delta$. 
As for the second condition notice that for any $\y\in \Xi$ 
we have $d(\y,\leps )\geq \epsilon -(2\alpha r+2\epsilon_0)\geq 3\alpha r$.
As $d(\Xi,\partial B(\x,r))=\frac{1}{6}r>2\alpha r$ 
the third condition is satisfied as well.
\end{proof}

\section{Friendly Schmidt's game}

Let $(X,d)$ be a complete metric space and let ${\cal{S}}\subset X$ be 
a given set (a target set). 
\textbf{Schmidt's game} \cite{S} is played by two players $A$ and $B$, 
each equipped with parameters $\alpha $ and $\beta $ 
respectively, $0<\alpha ,\beta <1$. 
The game starts with player $B$ choosing $y_0\in X$ and $r>0$ 
hence specifying a closed ball $B_0=B(y_0,r)$. 
Player $A$ may now choose any point $x_0\in X$ provided that 
$A_0=B(x_0,\alpha r)\subset B_0$. 
Next, player $B$ chooses a point $y_1\in X$ such that 
$B_1=B(y_1,(\alpha \beta)r)\subset A_0$. 
Continuing in the same manner we have a 
nested sequence of non-empty closed sets  
$B_0\supset A_0\supset B_1\supset A_1\supset ...\supset B_k\supset A_k...$ 
with diameters tending to zero as $k\rightarrow\infty$. 
As the game is played on a complete metric space, 
the intersection of these balls 
is a point $z\in X$. Call player $A$ the winner if $z\in {\cal{S}}$. 
Otherwise player $B$ is declared winner. 
A strategy consists of specifications for a player's choices 
of centers for his balls 
as a consequence of his opponent's previous moves. 
If for certain $\alpha $ and $\beta $ 
player A has a winning strategy, i.e., 
a strategy for winning the game regardless of how well player B plays, 
we say that ${\cal{S}}$ is an 
$\boldsymbol{(\alpha , \beta)}$\textbf{-winning} set. 
If it so happens that $\alpha $ is such that 
${\cal{S}}$ is an $(\alpha , \beta)$-winning set 
for all $0<\beta <1$, we say that ${\cal{S}}$ is an 
$\boldsymbol{\alpha }$\textbf{-winning} set. 
Call a set \textbf{winning} if such an $\alpha $ exists.

\vspace{1cm}

\noindent We define the following (target) set. 
This definition is a modification of the one given in \cite{KTV}.

\begin{defn}Suppose $\Omega \subset \RN$ and let 
${\cal U}= \{U_j \subset \mathbb{R}^N:j\in \mathbb{N}\}$  
be a family of subsets of $\mathbb{R}^N$. 
If $I:\mathbb{N} \rightarrow \mathbb{R}^+$ 
is an increasing function tending to infinity 
as $j$ tends to infinity and $\rho :\mathbb{R}^+\rightarrow \mathbb{R}^+$ 
is such that $\rho (r) \rightarrow 0$ as $r \rightarrow \infty$ 
and decreasing for large enough $r$, let

\begin{center}$\textbf{Bad} ^{*}({\cal {U}},I,\rho, 
\Omega)=\left\{\mathbf{x}\in \Omega :\exists \delta >0\text{ such that } 
d(\x,U_j) \geq \delta \rho (I(j))\; \forall j \in \mathbb{N}\right\}$.
\end{center}
\end{defn}    
As an immediate consequence of the above definition we get:
\begin{cor}\label{BA=BAD}
For $\Omega \subset\RN$, and $j\in\mathbb{N}^{+}$ defining 
$U_j=\left\{\frac{\p}{j}: \p\in \mathbb{Z}^{N}\right\}$, $I(j)=j$
 and $\rho (I(j))=j^{-{\frac{N+1}{N}}}$, we have
\begin{center}
$\textbf{BA}\cap \Omega=\textbf{Bad} ^{*}({\cal {U}},I,\rho,\Omega)$
\end{center}

\end{cor}

\noindent In the following theorem we shall show that under certain assumptions, 
$\textbf{Bad} ^*({\cal U},I,\rho, \Omega)$ is an $(\alpha , \beta)$-winning set.

\begin{theorem}\label{BAD=ALPHABETA} 
Suppose $\mu$ is absolutely friendly (with constants as in definition \ref{def:aad}) 
and $(\frac{D}{C})^{^{\frac{1}{a}}}= \alpha^{'}$. 
Let $\Omega=\K$ and suppose $F:\mathbb{N}\rightarrow \mathbb{R}^{+} $ 
is an increasing function, with $F(k)\rightarrow \infty $ as $k\rightarrow \infty $. 
Define $F^{0}=\left[0,F(0)\right)$ and $F^{k}=\left[F(k-1),F(k)\right)$ 
for any $k>0$. Let ${\cal U}= \{U_j \subset \mathbb{R}^N:j\in \mathbb{N}\}$  
be a family of subsets of $\mathbb{R}^N$. 

\noindent Suppose $0<\beta <1$ and $0<\alpha < \frac{1}{12}\alpha ^{'}$ satisfy:  
\begin{enumerate}

\item  for every $k,l\in \mathbb{N}$, for every $\x\in \K$ and for every $r\leq r_0$, 

if $I(j_1),...,I(j_l)\in F^{k}$ then 
$\left(\bigcup^{l}_{i=1}U_{j_{i}}\right)\bigcap B(\x ,(\alpha \beta)^{k}r)\subset {\cal L}$ 
for some affine hyperplane $\cal L$,

\item  for every $k$, $(\alpha \beta)^k\geq \rho (F(k))$. 

\end{enumerate}
Then $\textbf{Bad} ^*({\cal U},I,\rho,\Omega)$ is an 
$(\alpha , \beta )$-winning set on $\Omega $.
\end{theorem}

\begin{proof}
Player A's strategy is to play in an arbitrary manner until the the first ball of 
radius $r_I\leq r_0$ is chosen by player B. 
Let $k_0\in \mathbb{N}$ be such that 
$\beta ^{k_0+1}r_0<r_I\leq \beta ^{k_0}r_0$. 
Set $\delta=(\alpha\beta) ^{k_0+1}\beta^{k_0} r_0$ 
and let $r'=(\alpha\beta)^{k_0}r_I$.  
\\We ``reset'' our counter and specify player A's strategy from this point on. 
At his $k$th move player $A$ has to choose a point $\x\in \K$ such that 
$A_k=B(\x,\alpha (\alpha \beta)^{k}r')\subset B_k=B(\y,(\alpha \beta)^{k}r')$ 
where $\y\in\K$ is player B's $k$th choice. 
Let ${\cal{U}}_j=\bigcup^{l}_{i=1}U_{j_{i}}$ where $I(j_1),...,I(j_l)\in F^{k}$.
\begin{itemize}
\item[(a)]
If ${\cal{U}}_j\bigcap B(\y ,(\alpha \beta)^{k}r')=\emptyset$, 
player A may choose $\x=\y$.

By Lemma \ref{GEO}(3) 
\begin{center}
$d\left({\cal{U}}_j,A_k\right)>\alpha (\alpha\beta)^{k}r'\geq \delta 
(\alpha\beta)^{k}\geq \delta \rho (F(k))>\delta \rho (I(j)$.
\end{center}  
\item[(b)]
Otherwise suppose ${\cal{U}}_j\bigcap B(\y ,(\alpha \beta)^{k}r')\neq \emptyset$.

by Lemma \ref{GEO}(2) player A can pick a point $\x=\x_k$ such that 
\begin{center}
$d\left({\cal{U}}_j\bigcap B(\y ,(\alpha \beta)^{k}r'),A_k\right)>\alpha 
(\alpha\beta)^{k}r'>\delta\rho(I(i))$. 
\end{center}
Furthermore, if 
${\cal{U}}_j-B(\y ,(\alpha \beta)^{k}r')\neq \emptyset$ then by Lemma \ref{GEO}(3)

\begin{center}
$d\left({\cal{U}}_j-B\left(\y , (\alpha\beta)^{k}r'\right),A_k\right)>\alpha 
(\alpha\beta)^{k}r'>\delta\rho(I(i))$.
\end{center} 
\end{itemize}
\end{proof}

The following proposition due to W.M. Schmidt \cite{S} (Theorem 2) 
is material for later considerations.   
\begin{pro}\label{CIP}
The intersection of countably many $\alpha $-winning sets is $\alpha $-winning. 
\end{pro}

\section{Full Hausdorff dimension}

We now are in position to formulate a sufficient condition
for establishing a lower bound of a winning set's Hausdorff dimension,
where the winning set is a subset of 
the support of an absolutely friendly measure.

The main ideas in this section are due to W. M. Schmidt \cite{S}. 
We nonetheless have decided to include the definitions, results
and proofs for the sake of clearer understanding the 
connection to the previous definitions and results.

\begin{defn}For a metric space $(X,d)$, given $x\in X$, and real numbers $r>0$, $0<\beta <1$, 
denote by $N_{X}(\beta ,x,r )$ the maximum number of disjoint balls 
of radius $\beta r$ contained in $B(x,r)$.
\end{defn}

\newpage

\begin{theorem}\label{AF=LB} Let $\mu $ be absolutely friendly and denote $X=\K$. 
Suppose the following condition is satisfied:

There exists constants $r_1\leq 1$, $M$ and $\delta $ such that for every 
$0<r\leq r_1$, $0<\beta <1$ and $\x \in X $, 
\begin{equation}\label{no of balls}
N_{X}(\beta ,\x,r )  \geq M\beta ^{-\delta }. 
\end{equation}

Then if $\cal{S}$ is a winning set on $(X,d)$ then $\text{dim}\cal{S}\geq \delta $.  
\end{theorem}

In the course of the proof of we shall use the following auxiliary lemma.
(Lemma 20 in \cite{S}).

\begin{pro}\label{HS}
Let $\cal{H}$ be a Hilbert space and let $w_0=2\sqrt{3}-1$. 
For any $r\in \mathbb{R}^{+}$ let ${\cal{M}}$ be any collection of balls 
$\left\{B(x_i,r): i\in \mathbb{N},\; x_i\in {\cal{H}}\right\}$ such that 
\begin{center}
$\text{for every}\;\; i\neq j, \; \;int B(x_i,r)\cap int B(x_j,r)=\emptyset$.\\
\end{center}

Then for any $r_0<w_0r$ and $x\in {\cal{H}}$ the ball $B(x,r_0)$ 
has a non empty intersection with at most two balls from ${\cal{M}}$.
\end{pro}

\begin{proof}\textsl{{Theorem \ref{AF=LB}}}. 

Let $\mu $ be an absolutely friendly measure 
satisfying condition \ref{no of balls} and $\beta \leq (\frac{M}{2})^{\frac{1}{\delta }}$. 
Thus $\N\geq 2$ for every $\x\in \K $. 
In order to estimate the Hausdorff dimension of a winning set $\cal{S}$ 
assume player
$A$ is playing to win the game using some strategy.
This means that given choices of balls $B_0 \supset A_0 \supset \ldots
A_{k-1} \supset B_k$, played by the two players prior to player $A$'s kth turn, the
strategy of player $A$ chooses a ball $A_k \subset B_k$. Since the strategy is
winning, $\bigcap A_k = \bigcap B_k$ will be in $\cal{S}$ regardless of player $B$'s
choices. Here we will describe many possible strategies for player $B$, resulting
in many points in $\cal{S}$.

We consider the game from the loser's point of view, player $B$. 
Fix $\beta $ such that 

\begin{center}
$2\leq N(\beta )=\text{min}\left\{\N:\x\in X, \; 0<r\leq r_1 \right\}.$
\end{center}

At each stage of the game player $B$ may direct the game to $N(\beta )$ 
disjoint balls and we restrict his moves to these $N(\beta )$ choices. 
Thus for each sequence of choices made by player $B$ 
with the restriction above,
we obtain a parametrization of the sequence of balls chosen by him. 
Let $B_0$ be his initially chosen ball, and for $k\in \mathbb{N}^{+}$, 
corresponding to his kth move, let $B_k=B_k(j_1,...,j_k)$, 
with $j_i\in \left\{0,..., N(\beta ) -1\right\}$ $i=1,2,...,k$. 
Notice also that given a sequence of positive integers $i_1$,$i_2$,... 
there is a \textsl{unique} point $x=x(i_1,i_2,...)$ 
contained in \textsl{all} balls  $B_k=B_k(j_1,...,j_k)$. 
By considering the $N(\beta )$ ways in which player $B$ 
may direct the game we consider the function 

\begin{center}
$f:\left\{0,..., N(\beta ) -1\right\}^{\mathbb{N}}\rightarrow \cal{S}$, 
$(t _k)_{k\in \mathbb{N}}\mapsto \bigcap _{{{k\in \mathbb{N}}}}B_k(t _1,...,t_k)
=\left\{x(t )\right\}$.
\end{center} 
As every number in the closed unit interval has 
at least one expansion in base 
$N(\beta )$ we map the image of $f$, $\cal{S}^{*}\subset \cal{S}$ 
onto $\left[0,1\right]$ by 

\begin{center}
$g:{\cal{S}^{*}}\rightarrow [0,1]$, $x(t )\mapsto 0.t_1t_2 ...$.
\end{center}

In view of proposition \ref{HS}, for $0<w<w_0$ and $0<\alpha <1$ 
any ball of radius $w(\alpha \beta)^k$ 
intersects at most two of the balls $B_k(j_1,...,j_k)$. 
Let ${\cal C}=\left\{C_l\right\}_{l\in \mathbb{N}}$ 
be a cover of $\cal{S}\cap K$ 
of balls with radius $\rho (C_l)=\rho_l$. As $\cal C$ covers $\cal{S}^*$ 
we have that $g(\cal C)$ covers $\left[0,1\right]$. 
Let $ \overline{\lambda}$ denote the outer Lebesgue measure. 
We have
\begin{equation}\label{3.11}
\sum ^{\infty}_{l=1}\overline{\lambda}(g(C_l))
\geq \overline{\lambda}(\bigcup^{\infty}_{l=1}g(C_l))\geq 1.
\end{equation}
\noindent Define integers
\begin{center}
$k_l=\left[k^{*}_l\right]$ where 
$k^{*}_l=log_{\alpha \beta}(2w^{-1}\rho _l)$.
\end{center}

\noindent Notice that:\\\\ 

\noindent $(2w^{-1}\rho _l)^\frac{logN(\beta)}
{\left|log(\alpha \beta )\right|}=N(\beta )^{-k^{*}_l}$ 
and since $k^{*}_l<k_l+1$ we get
\begin{equation}
N(\beta ) ^{-k_l}<N(\beta )  N(\beta ) ^{-k^{*}_l}=N(\beta ) 
(2w^{-1}\rho _l)^\frac{logN(\beta)}
{\left|log(\alpha \beta )\right|}.
\end{equation}
 
\noindent Assuming without loss of generality that for every 
$l$,  $\rho _l\leq \frac{w}{2}$, 
there exists $n_0\in \mathbb{N}$ such that 
$\frac{w}{2}(\alpha \beta )^{n_0+1}<\rho_l\leq\frac{w}{2}(\alpha \beta )^{n_0}$. 
It follows that $k_l=n_0$ and so\\ 
\begin{equation}
\rho_l<w(\alpha \beta )^{k_l}.
\end{equation}

\noindent This implies that the ball $C_l$ 
intersects at most two of the balls 
$B_{l}(j_1,...,j_{k_l})$. As the length of the interval $g(B_l(j_1,...,j_l))$ 
is $N(\beta ) ^{-k_l}$ 
we have $\overline{\lambda }(g(C_l))\leq 2N(\beta ) ^{-k_l}$. 
Combining with \ref{3.11},\\

\noindent $1\leq \sum ^{\infty}_{l=1}\overline{\lambda}(g(C_l))\leq 
\sum ^{\infty}_{l=1}2N(\beta ) ^{-k_l}< 2N(\beta )  (2w^{-1})^{\frac {log(N(\beta ) )}
{\left|log(\alpha \beta)\right|}}\sum ^{\infty}_{l=1 }\rho_l^{\frac{log(N(\beta ) )}
{\left|log(\alpha \beta )\right|}}$.\\

\noindent By definition, dim${\cal{S}}\geq \frac{log(N(\beta ) )}
{\left|log(\alpha \beta )\right|}\geq \frac {\delta \left|logC_0\beta \right|}
{\left|log\alpha \right|+\left|log\beta \right|}\rightarrow \delta$ 
as $\beta \rightarrow 0$.
\end{proof}

\begin{remark}If it so happens that $\delta $=dim($\K $) then obviously 
\begin{center}
dim$\cal{S}$=$\delta $.
\end{center}
\end{remark}

\section{Simplex lemma}
Before giving our main example in the following section, 
we prove a version of the simplex lemma following 
ideas credited by W.M.Schmidt in \cite{S} to Davenport.
\begin{theorem}\label{SIMPLEX}

\noindent Let $\Lambda :\RN \rightarrow \RN $ 
be an affine map and denote by $\A$ 
the $N\times N$ matrix associated with the linear part of $\Lambda $. 
For every $\theta \in (0,1)$ let $R=\theta ^{\frac{-N}{N+1}}$ 
and for every $k\in \mathbb{N}^{+}$ let
\begin{center}
$U_k=\left\{\Lambda(\frac{\p}{q}):q\in \mathbb{N}^{+}, 
\p \in \mathbb{Z}^{N}\;\text{and}\;\; R^{k-1}\leq q<R^{k} \right\}$.
\end{center} 
Denote by $V_N$ the volume of the $N$-dimensional unit ball.
Then for every $r>0$ such that $r^N<\left|\text{det}\A\right|(N!)^{-1} V_N^{-1} \theta ^N$ 
and for every $\x $ there exists an affine hyperplane $\cal{}L$ such that 
\begin{center}$U_k\cap B(\x,\theta ^{k-1}r) \subset \cal{L}$.
\end{center}

\end{theorem}

\begin{proof}
Assume the contrary and let 
$\left\{V_i\right\}^{N}_{i=0}$, $V_i=(v^{1}_{i},...,v^{N}_{i})$ 
be $N+1$ independent points in $U_k\cap B(\x,\theta ^{k-1}r)$, i.e., 
not belonging to any single affine hyperplane. 
Denote by $\Delta $ the $N$-dimensional simplex subtended by them. 
By a well known result from calculus we have \\

\noindent $\lambda _N(\Delta )=(N!)^{-1}\left|\text{det}L^{'}\right|>0$, 
where $L^{'}=
\left(
\begin{matrix}
v^{1}_{1}-v^{1}_{0} & . & . & . & v^{N}_{1}-v^{N}_{0} \\
. & . & . & .& . \\
. & . & . & .& . \\
. & . & . & .& . \\
v^{1}_{N}-v^{1}_{0} & . & . & . & v^{N}_{N}-v^{N}_{0}
\end{matrix}\right)$.\\

\noindent As $\lambda _N(\Delta )>0$ we have $\text{det}L^{'}\neq 0$. 

\noindent Consider now the $(N+1\times N+1)$ matrix $L=
\left(
\begin{matrix}
1 & v^{1}_{0} & . & .& v^{N}_{0} \\
. & . & . & .& . \\
. & . & . & .& . \\
. & . & . & .& . \\
1 & v^{1}_{N} & . & .& v^{N}_{N}
\end{matrix}\right)$.

By repeatedly subtracting the first row from all others we get 
$\text{det}L=\text{det}L^{''}$ where $L^{''}=
\left(
\begin{matrix}
1 & v^{1}_{0} & . & .& v^{N}_{0} \\
0 & v^{1}_{1}-v^{1}_{0} & . & .& v^{N}_{1}-v^{N}_{0} \\
. & . & . & .& . \\
. & . & . & .& . \\
0 & v^{1}_{N}-v^{1}_{0} & . & .& v^{N}_{N}-v^{N}_{0}
\end{matrix}\right)$ 
and so $\text{det}L=\text{det}L^{'}$.
\\

\noindent Hence, 
\noindent $\lambda _N(\Delta )=\left|\text{det}A\right|(N!)^{-1}
\left|\text{det}L \right|$ 
\vspace{1cm}
where $L$=$\left(\begin{matrix}
1 & \frac{p_0^1}{q_0} & . & .& \frac{p_N^1}{q_0} \\
. & . & . & .& . \\
. & . & . & .& . \\
. & . & . & .& . \\
1 & \frac{p_0^N}{q_N} & . & .& \frac{p_N^N}{q_N}
\end{matrix}\right)$ 

and $\text{det}L\neq 0$ by our assumption.\\ 

\noindent Notice also that $q_0\cdot q_1\cdot ...\cdot q_N\cdot L= 
\left(\begin{matrix}
q_0 & p_0^{1} & . & .& p_N^{1} \\
. & . & . & .& . \\
. & . & . & .& . \\
. & . & . & .& . \\
q_N & p_0^{N} & . & .& p_N^{N}
\end{matrix}\right)$, 
\vspace{1cm}

and as all entries in $q_0\cdot q_1\cdot ...\cdot q_N\cdot L$ 
are integers it follows that  
\begin{center}
$q_0q_1\cdot ...\cdot q_N\cdot\left|\text{det} L\right|\geq 1$. 
\end{center}

\noindent And so, 
\begin{equation}\label{4.14}
\lambda _N(\Delta )=(N!)^{-1}\left|\text{det}A\right|\left|\text{det}L\right|
\geq (N!)^{-1}\frac{\left|\text{det}\A\right|}{q_0\cdot... \cdot q_N} >(N!)^{-1}\left|
\text{det}\A\right|R^{-k(N+1)}.
\end{equation}

\noindent But,
\begin{equation}
\lambda _N (B(\x,\theta ^{k-1}r))=
(\theta ^{k-1} r)^{N} V_N=
{\theta ^{(k-1) N}r^N} V_N<\left|\text{det}\A\right|\theta ^{k N}(N!)^{-1},
\end{equation} 
\begin{equation}
\theta ^{kN}=(\theta ^{\frac{-N}{N+1}})^{-k(N+1)}=R^{-k(N+1)},
\end{equation} 
and so
\begin{equation}\label{4.17}
\lambda _N (B(\x,\theta ^{k-1}r))\leq \left|\text{det}\A\right|(N!)^{-1}R^{-k(N+1)}.
\end{equation}

\noindent by our assumption on $U_k$. 

\vspace{1cm}

\noindent As $\Delta \subset B(\x,\theta ^{k-1}r)$, 
\ref{4.14} contradicts \ref{4.17}.
\end{proof}

\vspace{2cm}

\section{Application to Hutchinson's construction}

Before turning our attention to our main example
we state and prove the following theorem which is material
for what follows.

\begin{defn}
Say that $\mu $ satisfies the power law if 
there exist real numbers $a,b,\delta>0$ such for every 
$\x\in \K$, $0<r\leq 1$
\begin{center}
$a r^\delta \leq \mu (B(\x,r))\leq b r^\delta$.
\end{center}
\end{defn}

\begin{theorem}\label{PL implies 3.5}
Let $\mu $ satisfy the power law. 
Then $\mu$ 
satisfies condition \ref{no of balls}. 
\end{theorem}

\begin{proof}
\noindent Let $r\leq 1$, $0<\beta <1$ and consider a ball $B(\x,r)$ 
with $\x \in K$. Denote by $\left\{\x_i\right\}$, 
$i\in \left\{0,...,\N \right\}$ 
the centers of the $\N $ balls under consideration. 
Then, for every $i$, $\x_i\in B(\x,(1-\beta )r)\cap K$.

\noindent By a simple geometric argument we see that 
the collection of balls 
$B(\x_i,3\beta r)$ cover $B(\x,(1-\beta )r)$. 
For otherwise there exists $\y \in B(\x,(1-\beta )r)$ 
such that $d(\y,\x_i)\geq 3\beta r$ 
for every $i$. It follows that $B(\y,\beta r)$ could be added 
to the original collection of balls, 
which is a contradiction to the maximality assumption on $\N $. 
We may assume that $\beta\leq\frac{1}{2}$ 
with no loss of generality, 
as for $\frac{1}{2}<\beta<1$ we may choose 
$M\leq 2^{-\delta}\Rightarrow M\beta^{-\delta}\leq 1$. 
Notice also that $\delta \leq N$. And so,
\\\\
$a(1-\beta )^{\delta}r^{\delta }\leq \mu(B(\x,(1-\beta )r)
\leq \N \mu(B(\x_i,3\beta r))
\leq \N b 3^{\delta }\beta ^{\delta  }r^{\delta }$.

\begin{equation}
\N \geq  ab^{-1}3^{-1}(1-\beta )^{\delta }\beta ^{-\delta }
\geq ab^{-1}3^{-1}2^{-N}\beta ^{-\delta }.
\end{equation}
Thus condition \ref{no of balls} is satisfied with $r_1=1$ 
and $M=ab^{-1}3^{-1}2^{-N}$.

\end{proof}

\noindent A map $\phi : \RN \rightarrow \RN$ is a  \textbf{similarity}
if it can be written as
\begin{center}
$\phi (\x) = \rho \Theta (\x) + \y,$
\end{center}
where $\rho \in \mathbb R^{+}$, $\Theta \in O (N,\mathbb R)$ 
and $\y \in \RN $.
It is said to be \textbf {contracting} if $\rho < 1$.
It is known (see [Hu] for a more general statement) that for
any finite family
$\phi_1,\dots, \phi_m$ of  contracting similarities there exists a
unique nonempty compact set $K$,
called the \textbf{attractor} or \textbf{limit set} of the family,
such that
\begin{center}
$K = \bigcup_{ i = 1 }^m\phi_i ( K ).$
\end{center}
Say that
$\phi_1, \dots,  \phi_m$ as above satisfy the {\textbf{open set
condition} 
if there exists an open subset $U \subset \RN $ such that
\[\phi_i ( U ) \subset U \ \mathrm{for \ all \ \ } i=1, \ldots,
m\,,\]
          and
\[i \ne j \Longrightarrow \phi_i ( U ) \cap \phi_j(U) =
\varnothing\,.\]
The family $\{\phi_i\}$ is called \textbf{irreducible} if there is no
finite collection of proper affine subspaces which is invariant under
each $\phi_i$.
Well-known self-similar sets, like Cantor's ternary set, Koch's curve
or Sierpinski's gasket, are all examples of attractors of irreducible
families of
contracting similarities satisfying the open set condition.

Suppose $\{\phi_i\}_{i = 1}^m$ is a family of contracting 
similarities of  $\RN$
satisfying the open set
condition, let $K$ be its attractor, $\delta $ 
the Hausdorff dimension of
$K$, and $\mu$ the restriction of the $\delta $-dimensional Hausdorff
measure to $K$. 

J.~Hutchinson \cite{H} gave a simple formula for calculating
$\delta $ and proved that $\mu(K)$ is
positive and finite. Furthermore, 

\begin{pro}\label{PL} 
$\mu$ satisfies 
the power law with $\delta$=dimK.
\end{pro}

\noindent As a consequence of proposition \ref{PL} and 
theorem \ref{PL implies 3.5} we prove the following.
 
\begin{cor}\label{4.3}
Let $\left\{ \phi_1,...,\phi_k\right)$ be a finite irreducible family of 
contracting similarity maps of $\RN $ satisfying the 
open set condition. Let $K$ be its attractor. 
Let $\mu $ be the restriction of $H^{\delta }$ to $K$.
Then $\mu $ is absolutely friendly satisfying 
condition \ref{no of balls} with dim$K=\delta $.
\end{cor}

\begin{proof}
By theorem \ref{PL implies 3.5}, condition \ref{no of balls} is satisfied.

Set $r_0=1$. 
It is easily seen that the power law implies that condition $(ii)$
of definition \ref{def:aad} is satisfied
with $D=\frac{a}{b} (\frac{5}{6})^{\delta }$.

Following  \cite{KLW}(Theorem 2.3, Lemma 8.2 and 8.3), there exist 
$C$ and $a$ such that $\mu$ is absolutely $(C,a)$-decaying (see remark 2) on 
any ball of radius $r=1$ centered in $\K$. 

Using the notation of Definition \ref{def:aad}, $\mu $ 
is absolutely friendly with $r_0=1$.

\end{proof}

We are now ready to prove our main example.

\begin{cor}\label{MAIN EXAMPLE}
Let $\left\{ \phi_1,...,\phi_k\right)$ be a finite irreducible family of 
contracting similarity maps of $\RN $ satisfying the open set condition. 
Let $K$ be its attractor and $\alpha^{'}$ as in lemma \ref{GEO}.
\noindent Then for any countable collection of affine transformations 
\noindent $\left\{\Lambda_i\right\}^{\infty}_{i=0}$, 
with $\Lambda_i:\RN \rightarrow \RN $ the set
\begin{center} 
${\cal{S}}= K\cap (\cap ^{\infty}_{i=1}( \Lambda_i(\textbf{BA})))$ 
\end{center} 
is an $\alpha $-winning set on $K$ for any 
$0<\alpha < \frac{1}{12}\alpha ^{'}$. 
Furthermore, dim${\cal \cal{S}}$=dim$K$.
\end{cor}
\begin{proof}In view of proposition \ref{CIP} it suffices to prove 
that for each $i$, $K\cap \Lambda_i(\textbf{BA})$ 
is $\alpha $-winning. 
Given an affine transformation $\Lambda $ and following 
corollary \ref{BA=BAD} we prove that 
$\textbf{Bad} ^{*}({\cal {U}},I,\rho, \Omega)$ 
is an $\alpha $ winning set on $\Omega =K$ 
where for every $q\in \mathbb{N}^{+}$
\begin{equation}
U_q=\left\{\Lambda (\frac{\p}{q}):\p\in \mathbb{Z}^{N}\right\},
\end{equation}
$I(q)=q$ and $\rho (I(q))=\rho(q)=q^{\frac{-N+1}{N}}$.
Following the notation of theorem \ref{BAD=ALPHABETA} and theorem \ref{SIMPLEX} 
let $\theta =\alpha \beta$ and for every $k\in \mathbb{N}^{+}$ 
let $F(k)=R^{k}=(\alpha \beta)^{\frac{-Nk}{N+1}}$.
Define
\begin{equation}
U_k=\left\{\Lambda (\frac{\p}{q}):q\in \mathbb{N}^{+}, 
\p\in \mathbb{Z}^{N}\; \text{and}\; R^{k-1}\leq q <R^{k}\right\}.
\end{equation}

By Theorem \ref{SIMPLEX} we get that the first condition of 
theorem \ref{BAD=ALPHABETA} 
is satisfied by any $\beta $.
As by our definition $\rho(F(k))=(\alpha \beta)^{k}$, 
the second condition is satisfied as well. 
Thus $K\cap T_i(\textbf{BA})$ is an $(\alpha, \beta)$-winning set for every $\beta $, 
rendering it an $\alpha $-winning set. 

Furthermore, as $\mu $ is absolutely friendly satisfying 
condition \ref{no of balls} with the exponent of the condition 
being $\delta $=dim$K$, by theorem \ref{AF=LB}, 
followed by remark 3 we are done. 

\end{proof}

\bibliographystyle{alpha}

\end{document}